\documentclass[11pt,reqno]{amsart}

\usepackage{xcolor}
\usepackage[utf8]{inputenc}
\usepackage[T1]{fontenc}
\usepackage[english]{babel}

\usepackage{microtype}
\usepackage{amsmath,amssymb,amsthm}
\usepackage{amsthm, thmtools, thm-restate}
\usepackage{mathtools} 
\usepackage{bm}        
\usepackage{enumitem}

\usepackage{graphicx}  
\usepackage{caption}
\usepackage{subcaption}
\usepackage{tikz}      
\usepackage{mathrsfs}  
\usepackage{booktabs}  
\usepackage{geometry}
\usepackage{cancel}

\theoremstyle{plain}
\newtheorem{theorem}{Theorem}[section]      
\newtheorem{lemma}[theorem]{Lemma}          
\newtheorem{proposition}[theorem]{Proposition}
\newtheorem{corollary}[theorem]{Corollary}

\newtheorem*{theorem*}{Theorem}             

\theoremstyle{definition}
\newtheorem{definition}[theorem]{Definition}

\newtheorem{remark}[theorem]{Remark}
\usepackage[hypertexnames=false]{hyperref} 


\numberwithin{equation}{section}

\newcommand{\R}{\mathbb{R}}

\usepackage{comment}
\usepackage{graphicx}  
\usepackage{caption}
\usepackage{subcaption}
\usepackage{tikz}      
\usepackage{mathrsfs}  
\usepackage{booktabs}  
\usepackage{geometry}
\usepackage{cancel}

\title[Multiplicity for the Minkowski mean curvature equation]{Multiplicity via symmetry for the Minkowski mean curvature equation: Dirichlet and Neumann problems}

\author[A. Boscaggin]{Alberto Boscaggin}

\author[R. Ziegele]{Ricardo Ziegele}
\address{Alberto Boscaggin, Ricardo Ziegele \newline \indent 
Dipartimento di Matematica
\newline\indent
Universit\`a di Torino
\newline\indent
via Carlo Alberto 10, 10123 Torino, Italy}
\email{alberto.boscaggin@unito.it}
\email{ricardoalfonso.ziegelealiaga@unito.it}

\begin{document}
\thanks{Work written under the auspices of INdAM - GNAMPA}

\begin{abstract}
 We consider the Minkowski mean curvature equation
$$ -\operatorname{div}\left(\frac{\nabla u}{\sqrt{1-|\nabla u|^2}}\right)
=\lambda u+\mu h(x,u) \quad \text{in } \Omega, \quad \mathcal B u = 0 \quad \text{on } \partial \Omega,
$$
where $\mathcal B \in \{ \mathcal D, \mathcal N \}$ is a Dirichlet or Neumann boundary condition and $\Omega$ is a smooth bounded domain, required to be convex in the Neumann case.
Moreover, $\lambda, \mu \geq0$ are parameters, and $h(x,u)$ is odd with respect to $u$ and positive for $u > 0$. Using minmax methods in nonsmooth critical point theory, we prove that $\lambda>\lambda^\mathcal B_k$ yields at least $k$ pairs of nontrivial solutions with negative energy for every $\mu\geq0$, where $\lambda_k^\mathcal B$ are the corresponding Laplacian eigenvalues, indexed by $\lambda^\mathcal B_0=0$. Moreover, if $H(x,u):=\int_0^u h(x,t)\,dt=o(u^2)$ as $u\to0$, we prove that for $\lambda\in[\lambda^\mathcal B_k,\lambda^\mathcal B_{k+1})$ and every $j\geq\max\{1,k\}$, sufficiently large $\mu$ yields at least $2j-k$ distinct pairs of solutions: $j$ with negative energy and $j-k$ with positive energy.
\end{abstract}
\maketitle

\section{Introduction and statement of the main results}
In this paper, we are concerned with \emph{multiplicity results} for the Minkowski-curvature boundary value problem
\begin{equation}
\label{eq:problem}
\begin{cases}
-\displaystyle\operatorname{div}
\left(\dfrac{\nabla u}{\sqrt{1-|\nabla u|^2}}\right)
=\lambda u+\mu h(x,u)
& \text{in }\Omega,\\[6pt]
\mathcal B u=0
& \text{on }\partial\Omega,
\end{cases}
\end{equation}
where $\mathcal B\in\{\mathcal D,\mathcal N\}$ stands for the Dirichlet
or Neumann boundary operator, respectively, and
\begin{itemize}
    \item[$(H_\mathcal B)$] $\Omega \subset \R^N$ is a bounded $C^2$ domain, additionally required to be convex if $\mathcal{B} = \mathcal{N}$.
\end{itemize}
Finally, $\lambda,\mu\geq0$ are parameters and $h: \Omega \times \mathbb{R} \to \mathbb{R}$ is a function which satisfies:
\begin{itemize}
\item[$(h_1)$] $h$ is Carath\'eodory and for every $\varrho>0$ there exists $\alpha_\varrho \in L^\infty(\Omega)$ such that
$$
\ |h(x,s)| \le \alpha_\varrho(x)
\quad
\text{for a.e. } x\in\Omega\text{ and for all } |s|\le \varrho;
$$
\item[$(h_2)$] for a.e. $x \in \Omega$ and every $s > 0$, it holds that
$$0 < h(x,s) = -h(x,-s);$$
\end{itemize}
The right-hand side of the equation in \eqref{eq:problem} is thus an \emph{odd} function with respect to $u$. 

Existence and multiplicity questions for the Minkowski
mean curvature equation have been investigated by
several methods, particularly under Dirichlet boundary
conditions. Among these, we mention topological degree
arguments \cite{CoObOmRi13b,Da16} and variational
techniques \cite{BeJeMa14,BeJeTo13, GaLo26}. In the latter approach,
the gradient constraint $\Vert \nabla u \Vert_{L^\infty} \leq 1$
naturally leads to a nonsmooth action functional,
which can be studied within Szulkin's critical point
theory \cite{Sz86}. By contrast, far fewer results are available for 
Neumann problems, owing in
particular to a regularity issue that has only recently
been resolved for convex domains \cite{Ma26}, as we shall explain in more detail
below.

When the nonlinear term is odd, the variational approach
offers an additional mechanism for multiplicity:
the associated functional is even, and the symmetry
\(u \mapsto -u\) can be exploited through symmetric
minmax methods to obtain multiple pairs of nontrivial
solutions. This idea was used by Jebelean and \c{S}erban
\cite{JS} in the Dirichlet setting, under suitable conditions involving
the eigenvalues of the Laplacian. Along this line of research, in the present paper we establish new multiplicity
results for both Dirichlet and Neumann problems.

In order to give the statements of our results, let us denote by
\[
0<\lambda_1^{\mathcal D} <\lambda_2^{\mathcal D}
\leq
\lambda_3^{\mathcal D}
\leq \cdots\leq
\lambda_i^{\mathcal D}
\leq\cdots
\]
and by 
\[
0=\lambda_0^{\mathcal N} <\lambda_1^{\mathcal N}
\leq
\lambda_2^{\mathcal N}
\leq \cdots\leq
\lambda_i^{\mathcal N}
\leq\cdots
\]
the eigenvalues of the linear operator $-\Delta$ under homogeneous Dirichlet
and Neumann boundary conditions, respectively. Moreover, we set $\lambda_0^{\mathcal{D}} := 0$, so that in any case we have a sequence 
$\{\lambda_i^{\mathcal{B}}\}_{i = 0}^\infty$ with $\lambda_0^{\mathcal{B}} = 0$, both when $\mathcal{B} = \mathcal{D}$ and when $\mathcal{B} = \mathcal{N}$ (notice, however, that $\lambda_0^{\mathcal{D}} = 0$ by definition, while $\lambda_0^{\mathcal{N}} = 0$ is actually the first eigenvalue with Neumann boundary conditions).
Finally, throughout the paper the energy of a solution of \eqref{eq:problem} 
is understood as the value of the associated
variational functional, normalized to vanish
at the trivial solution (see Section \ref{sec-variational}).

\medbreak
Our first result reads as follows.

\begin{theorem}\label{th-spectral}
Assume $(H_\mathcal B)$ with  $\mathcal B\in\{\mathcal D,\mathcal N\}$, together with $(h_1)$ and $(h_2)$. Moreover, suppose that
\begin{equation}\label{hp-spectral}
\lambda > \lambda_k^{\mathcal{B}}, \quad \text{ for some $k \geq 1$.}
\end{equation}
Then for every $\mu \geq 0$
problem \eqref{eq:problem} admits at least $k$ distinct
    pairs of nontrivial solutions with negative energy.
    Moreover, all these solutions are sign-changing if $\mathcal{B} = \mathcal{N}$, while for $\mathcal{B} = \mathcal{D}$ at least a pair of constant-sign solutions exists. 
\end{theorem}

\begin{remark}
A couple of comments about the statement of Theorem \ref{th-spectral} are in order. 
First, we emphasize that assumption \eqref{hp-spectral} does not require $\lambda_k^\mathcal{B}$ to be a simple eigenvalue. Of course, in case of a multiple eigenvalue 
$\lambda_k^\mathcal{B} = \lambda_{k+\ell}^\mathcal{B}$ for some $\ell \geq 1$, at least $k + \ell$ pairs of nontrivial solutions can be obtained. That is to say, Theorem \ref{th-spectral} actually provides at least $k(\lambda)$ pairs of solutions, where 
\[
k(\lambda) := \#\{ k \geq 1 \, : \, \lambda_k^\mathcal{B} < \lambda\}.
\]
Second, we observe that, while in the Dirichlet case the existence of solutions is guaranteed as soon as $\lambda$ is larger than the first eigenvalue $\lambda_{1}^\mathcal{D}$, the existence of Neumann solutions is obtained when $\lambda$ is greater than the \emph{second} eigenvalue $\lambda_{1}^\mathcal{N}$. This can be explained by the fact that Neumann solutions must be sign-changing, while the Dirichlet solutions coming from the spectral conditions $\lambda > \lambda_{1}^\mathcal{D}$ are the ones which are of constant sign. 
The convention used to count the eigenvalues, however, allows us to give the statement of Theorem \ref{th-spectral} in a completely symmetric way for both $\mathcal{B} = \mathcal{D}$ and $\mathcal{B} = \mathcal{N}$.
\end{remark}

It is worth noticing that Theorem \ref{th-spectral} can be applied, in particular, when $\mu = 0$, therefore giving the following direct corollary. 

\begin{corollary}\label{cor-spectral}
Assume $(H_\mathcal B)$ with  $\mathcal B\in\{\mathcal D,\mathcal N\}$ and that \eqref{hp-spectral} holds true.
Then the problem \begin{equation}
\label{eq:problem_spectral}
\begin{cases}
-\displaystyle\operatorname{div}
\left(\dfrac{\nabla u}{\sqrt{1-|\nabla u|^2}}\right)
=\lambda u
& \text{in }\Omega,\\[6pt]
\mathcal B u=0
& \text{on }\partial\Omega,
\end{cases}
\end{equation}
admits at least $k$ distinct
    pairs of nontrivial solutions with negative
    energy. Moreover, all these solutions are sign-changing if $\mathcal{B} = \mathcal{N}$, while for $\mathcal{B} = \mathcal{D}$ at least a pair of constant-sign solutions exists. 
\end{corollary}

The purely spectral problem \eqref{eq:problem_spectral} in the Dirichlet case $\mathcal{B} = \mathcal{D}$ was first considered in the paper \cite{Mi21}, providing the existence of a pair of constant-sign solutions for any $\lambda > \lambda_1^{\mathcal{D}}$. The existence of $k \geq 2$ pairs of solutions was later obtained in \cite{JS} under the condition $\lambda>2\lambda_k^{\mathcal D}$. It is worth mentioning that the possibility of establishing such a result under the natural condition $\lambda>\lambda_k^{\mathcal D}$ was explicitly raised as an open question in \cite[Remark 2.6]{JS}. Corollary \ref{cor-spectral} shows that such a condition is indeed
sufficient (at least in the setting of problem \eqref{eq:problem_spectral}), for both Dirichlet and Neumann boundary conditions.

\begin{remark}
We observe that, when dealing with problem \eqref{eq:problem_spectral}, some of the spectral thresholds provided by Corollary \ref{cor-spectral} are also necessary. More precisely:
\begin{itemize}
\item for the Dirichlet problem, nontrivial solutions exist if and only if \(\lambda>\lambda_1^\mathcal{D}\), while sign-changing solutions exist if and only if \(\lambda>\lambda_2^\mathcal{D}\). Indeed, the necessity of the first condition follows by testing the equation with \(u\), whereas the second follows by testing separately with \(u^+\) and \(u^-\) and using the minmax characterization of \(\lambda_2^\mathcal{D}\). Conversely, when \(\lambda>\lambda_2^\mathcal{D}\), Corollary \ref{cor-spectral} provides at least two distinct pairs of solutions, and the uniqueness of the positive solution established in \cite{BoFoMoSaTa18} ensures that at least one of these pairs is sign-changing. In particular, for \(\lambda_1^\mathcal{D}<\lambda\leq\lambda_2^\mathcal{D}\), the problem admits exactly one pair of nontrivial solutions, consisting of a positive solution and its opposite.
\item for the Neumann problem, every nonconstant solution must change sign and, for \(\lambda>0\), have zero mean. Testing the equation with \(u\) and using the Poincaré–Wirtinger inequality shows that necessarily \(\lambda>\lambda_1^\mathcal{N}\). Since the converse follows from Corollary \ref{cor-spectral}, the existence of nonconstant solutions is equivalent to this condition. When \(\lambda=0\), the only solutions are the constant functions.
\end{itemize}
\end{remark}

To state our second main result, we introduce the following additional condition concerning the behavior of the nonlinear term near zero:
\begin{itemize}
\item[$(h_3)$] it holds that 
\[\lim_{s\to 0}\frac{H(x,s)} {s^2} = 0, \qquad \text{ uniformly for a.e. $x\in\Omega$},
\]
where $H(x,s) := \int_0^s h(x,\tau)\,d\tau$.
\end{itemize}
Under this assumption, a different
multiplicity result can be obtained for sufficiently
large values of $\mu$, involving solutions at both
negative and positive energy levels.
\begin{theorem}\label{th-mularge}
Assume $(H_\mathcal B)$ with  $\mathcal B\in\{\mathcal D,\mathcal N\}$, together with $(h_1)$, $(h_2)$ and $(h_3)$. Moreover, suppose that
\[
\lambda\in
    [\lambda_k^{\mathcal B},\lambda_{k+1}^{\mathcal B}), \quad \text{ for some $k \geq 0$}.
\]
Then for every integer $j\geq\max\{1,k\}$ there exists
    $\mu_{j,\mathcal B}^\star(\lambda)>0$ such that,
    for every $\mu>\mu_{j,\mathcal B}^\star(\lambda)$,
    problem \eqref{eq:problem} admits at least
    $2j-k$
    distinct pairs of nontrivial solutions, $j$ of which with negative energy and (if $j>k$) $j-k$ with positive energy.
     Moreover, all these solutions are sign-changing
if $\mathcal B=\mathcal N$, while for
$\mathcal B=\mathcal D$ at least one pair
of negative-energy solutions has constant sign.
\end{theorem}

The multiplicity count in Theorem \ref{th-mularge} parallels the one obtained in the radial setting in \cite{BoCoZi26} with a Dirichlet boundary condition. There, a shooting argument yields $j$ \emph{large} solutions and $j-k$ \emph{small} solutions, distinguished by their initial values and classified according to their number of zeros. Here, a variational approach yields the same total count of $2j-k$ pairs, with an analogous splitting into $j$ pairs at negative energy and $j-k$ pairs at positive energy. It is therefore tempting to interpret the sign of the energy as a variational counterpart of the distinction between large and small solutions in the radial setting, associating negative energy with the former and positive energy with the latter. This interpretation is supported, on the positive-energy side, by \cite[Proposition~2.17]{BoCoZi26}: in the Dirichlet setting and under the hypotheses of that result, any family of solutions whose energies remain bounded from below converges uniformly to zero as $\mu\to+\infty$. In particular, this applies to these positive-energy solutions and provides an asymptotic justification for viewing them as small solutions. The agreement between the multiplicity counts thus suggests a connection between the shooting and variational descriptions, although identifying the corresponding families, especially on the negative-energy side, requires further analysis.

Our proofs are based on two symmetric minmax
principles for nonsmooth functionals, developed
by Arcoya, Bereanu and Torres \cite{ABT}.
The first, relying on the Krasnosel'skii genus,
provides critical points at negative energy
levels, while the second exploits a symmetric
linking geometry to obtain additional critical
points at positive energy levels.
The key point in applying these abstract results is to establish
the corresponding minimax geometries through spectral estimates
on suitable finite-dimensional spaces and their orthogonal complements.
The variational framework is completed by the verification of a
weak Palais--Smale condition and by the use of existing regularity
results identifying critical points with strictly spacelike solutions.
In particular, the recent result of Maniscalco \cite{Ma26}
allows us to treat Neumann boundary conditions on convex domains.

The plan of the paper is the following.
In Section \ref{sec-abstract}, we recall the abstract setting and the two minmax results from \cite{ABT}.
In Section \ref{sec-variational}, we introduce the nonsmooth variational framework for problem \eqref{eq:problem}, establish the correspondence between critical points and solutions, and verify the weak Palais--Smale condition in this specific setting.
Finally, Sections \ref{sec-3} and \ref{sec-4} are devoted to the proofs of Theorems \ref{th-spectral} and \ref{th-mularge}, respectively.

\medbreak
\noindent
\textbf{Acknowledgments.} We thank Francesca Colasuonno for useful discussions.

\medbreak
\noindent
\textbf{Declaration on the use of AI tools.}
 ChatGPT 6 (OpenAI) was used as an auxiliary tool to assist in checking the presentation and internal consistency of the proofs, and in the linguistic and stylistic revision of the manuscript. All mathematical arguments and all AI-assisted suggestions were independently checked and validated by the authors, who take full responsibility for the content of the paper.

\section{Preliminaries}

\subsection{Abstract setting} \label{sec-abstract}
We work within the abstract framework developed in \cite{ABT}. Precisely, let $X$ be a real Banach space and 
\[
I=\Psi+\Phi:X\longrightarrow(-\infty,+\infty],
\]
be a functional such that:
\begin{itemize}
\item[(i)] $\Psi:X\to(-\infty,+\infty]$ is convex, with non-empty and closed effective
domain
\[
D(\Psi):=\{u\in X:\Psi(u)<+\infty\},
\] and the restriction of $\Psi$ to $D(\Psi)$ is continuous,
\item[(ii)] $\Phi \in C^1(X,\mathbb{R})$,
\item[(iii)] $\Psi$ and $\Phi$ are \emph{even} and
$I(0)=0$.
\end{itemize}

\begin{remark}\label{rem-functional-abstract}
We recall that, in Szulkin's nonsmooth critical point theory \cite{Sz86}, the convex functional \(\Psi\) is required to be lower semicontinuous on \(X\): in the above setting, this property follows from assumption (i), since \(D(\Psi)\) is closed and the restriction of \(\Psi\) to \(D(\Psi)\) is continuous. The stronger continuity assumption, required by the abstract results in \cite{ABT}, will play a role in the choice of the functional setting for our applications (see Remark \ref{rem-functional-application}).
\end{remark}

The notion of critical point of the functional $I$ is the typical one in Szulkin's theory and it reads as follows.

\begin{definition}\label{def-critical}
    A point $u\in D(\Psi)$ is called a \emph{critical point} of $I$
if
\[
\Psi(v)-\Psi(u)+\Phi'(u)[v-u]\geq0
\qquad\text{for every }v\in D(\Psi).
\]
\end{definition}

In contrast, instead of the usual Palais-Smale condition, a notion of weak Palais-Smale condition is formulated with respect to a real Banach space $Y$ into which $X$ is
continuously embedded. Let us stress that the topologies of both spaces play a role: the error in the variational
inequality is measured using the norm of $X$, whereas convergence
is required to be with respect to the norm of $Y$.

\begin{definition}\label{def:weak-PS}
$I$ satisfies the \emph{weak Palais--Smale
condition with respect to $Y$}, abbreviated as
$\mathrm{(wPS)}$, if every sequence
$(u_n)\subset D(\Psi)$ such that $I(u_n)\longrightarrow c\in\mathbb R$ and, for some sequence $\varepsilon_n\downarrow 0$,
\[
\Psi(v)-\Psi(u_n)+\Phi'(u_n)[v-u_n]
\geq -\varepsilon_n\|v-u_n\|_X
\qquad\text{for every }v\in D(\Psi),
\]
admits a subsequence converging in $Y$ to a critical point
$u\in D(\Psi)$ satisfying $I(u)=c$.
\end{definition}

We now state the two complementary multiplicity results 
for even nonsmooth functionals that will be used for the proof of Theorems \ref{th-spectral} and \ref{th-mularge}. 
The first is based on the Krasnosel'skii genus,
while the second is a symmetric linking theorem
of generalized mountain-pass type.
Both results are formulated under the weak
Palais--Smale condition introduced above.
Throughout the following statements,
assumptions (i)--(iii) are understood.

\begin{theorem}[{\cite[Corollary~1]{ABT}}]
\label{thm:ABT-negative}
Assume that $I$ satisfies $\mathrm{(wPS)}$, and let
$1\leq\ell\leq m$ be integers.
Suppose that there exists a closed linear subspace
$Z\subset X$ and a symmetric set $A\subset D(\Psi)$ such that
\[
\operatorname{codim}_X Z=\ell-1,
\qquad
\inf_{u\in Z}I(u)>-\infty,
\]
and
\[
A\ \text{is oddly homeomorphic to } \mathbb S^{m-1},
\qquad
\sup_{u\in A}I(u)<0.
\]
Then $I$ has at least $m-\ell+1$ distinct pairs of nonzero critical points with negative
critical values. More precisely, there exist minmax levels
\[
-\infty < c_\ell \leq \cdots \leq c_m < 0,
\]
each of which is a critical value of \(I\).
Moreover, if \(c_i=c_j=c\) for some
\(\ell\leq i<j\leq m\), then \(I\) has infinitely
many distinct pairs of critical points at the level \(c\).
\end{theorem}

\begin{theorem}[{\cite[Theorem~2]{ABT}}]
\label{thm:ABT-positive}
Assume that $I$ satisfies $\mathrm{(wPS)}$, and let
$0\leq\ell<m$ be integers.
Suppose that there exist linear subspaces $E,Z\subset X$,
with $Z$ closed, and constants $r>\rho>0$ and $\alpha>0$
such that
\[
\dim E=m,
\qquad
\sup_{\substack{u\in E\\ \|u\|_Y=r}}I(u)<0,
\]
and
\[
\operatorname{codim}_X Z=\ell,
\qquad
\inf_{\substack{u\in Z\\ \|u\|_Y=\rho}}I(u)\geq\alpha.
\]
Then $I$ has at least $m-\ell$ distinct pairs of nonzero critical points with critical
values greater than or equal to $\alpha$.
More precisely, there exist minmax levels
\[
0<\alpha\leq b_{\ell+1}\leq\cdots\leq b_m<+\infty,
\]
each of which is a critical value of \(I\).
Moreover, if \(b_i=b_j=b\) for some
\(\ell+1\leq i<j\leq m\), then \(I\) has infinitely
many distinct pairs of critical points at the level \(b\).
\end{theorem}

The critical values appearing in the preceding results
admit minmax characterizations, given respectively
by the classical Krasnosel'skii genus and by the
symmetric linking families introduced in
\cite[Section 2.2]{ABT}. We refer to that paper for the
precise definitions. We only point out that the first
genus level satisfies \(c_1=\inf_X I\), possibly
equal to \(-\infty\).

\subsection{Variational framework} \label{sec-variational}

To apply the abstract results of Section~\ref{sec-abstract},
we set
\[
X:=W^{1,\infty}(\Omega),
\qquad
Y:=C(\overline\Omega),
\qquad
\|u\|_X
:=
\|u\|_{L^\infty(\Omega)}
+\|\nabla u\|_{L^\infty(\Omega)}.
\]
Moreover, 
for $\mathcal B\in\{\mathcal D,\mathcal N\}$, we define
\[
\Psi_{\mathcal B}(u)
:=
\begin{cases}
\displaystyle
\int_\Omega
\left(1-\sqrt{1-|\nabla u|^2}\right)\,dx,
& u\in K_{\mathcal B},\\[1ex]
+\infty,
& u\in X\setminus K_{\mathcal B},
\end{cases}
\]
where
\[
K_{\mathcal D}
:=
\left\{
u\in X:
\|\nabla u\|_{L^\infty(\Omega)}\leq1,\
u|_{\partial\Omega}=0
\right\},
\qquad
K_{\mathcal N}
:=
\left\{
u\in X:
\|\nabla u\|_{L^\infty(\Omega)}\leq1
\right\}.
\]
Note that, in the Dirichlet case, the boundary conditions are incorporated into the effective domain of the convex
functional $\Psi$. As it is easy to check, $\Psi_{\mathcal B}$ is proper and convex, its effective domain $K_{\mathcal B}$ is
closed in $X$, and the restriction of $\Psi$ to $K_{\mathcal B}$
is continuous in the norm of $X$.
Hence, assumption (i) of Section \ref{sec-abstract} is satisfied.

As for the smooth part, we first define,
\[
\widetilde h(x,s)
:=
\begin{cases}
h(x,-d), & s<-d,\\
h(x,s),  & |s|\leq d,\\
h(x,d),  & s>d
\end{cases}
\qquad \text{ where } \, d := \operatorname{diam}(\Omega),
\]
and
\[
\widetilde H(x,s):=\int_0^s \widetilde h(x,t)\,dt.
\]
Notice that \(\widetilde h\) still satisfies
assumption \((h_2)\), so that
\(\widetilde H(x,\cdot)\) is even and nonnegative. We then define
\[
\Phi_{\lambda,\mu}(u)
:=
-\frac{\lambda}{2}\int_\Omega u^2\,dx
-\mu\int_\Omega \widetilde H(x,u)\,dx.
\]
Of course, $\Phi_{\lambda,\mu}\in C^1 (X,\mathbb{R})$, with 
\[
(\Phi_{\lambda,\mu})'(u)[w]
=
-\int_\Omega
\bigl(\lambda u+\mu  \widetilde h(x,u)\bigr)w\,dx,
\]
so that assumption (ii) of Section \ref{sec-abstract} is satisfied as well.
Finally, both $\Psi_{\mathcal B}$ and $\Phi_{\lambda,\mu}$
are even, and, for $I_{\lambda,\mu}^{\mathcal B}
:=
\Psi_{\mathcal B}+\Phi_{\lambda,\mu}$, we have
\[
I_{\lambda,\mu}^{\mathcal B}(0)=0,
\]
so that assumption (iii) holds true.
Thus, all the structural assumptions of
Section~\ref{sec-abstract} are satisfied.
The weak Palais--Smale condition will be proved to hold in Lemma \ref{lem:compactness} below.

\begin{remark}\label{rem-functional-application}
It is worth observing that in previous variational treatments of the Minkowski mean curvature operator (see, e.g., \cite{BeJeMa14, BoCoZi26}) the energy functional was defined on \(C(\overline\Omega)\), with its convex part taking the value \(+\infty\) outside the admissible set \(K_{\mathcal B}\). This setting is suitable for Szulkin's critical point theory, since it can be proved that \(\Psi_B\) is lower semicontinuous with respect to uniform convergence. 
However, the abstract multiplicity results in \cite{ABT} recalled in Section \ref{sec-abstract} require the additional property that \(\Psi_B\) must be continuous on its effective domain. This is not true in the topology of \(C(\overline\Omega)\), since uniform convergence of admissible functions does not generally imply convergence of their gradients. On the other hand, such continuity holds in the stronger topology of \(W^{1,\infty}(\Omega)\), which motivates our choice of \(X\).
At the same time, compactness is 
naturally available in \(C(\overline\Omega)\), rather than in \(W^{1,\infty}(\Omega)\). This explains why we work with the continuous embedding \(X\hookrightarrow Y=C(\overline\Omega)\) and verify the weak Palais–Smale condition with respect to \(Y\), instead of the usual Palais–Smale condition in \(X\).
\end{remark}

Critical points \(u\in K_B\) of the functional
\(I_{\lambda,\mu}^B\) are understood in the sense
of Definition \ref{def-critical}, that is,
\[
\Psi_B(v)-\Psi_B(u)
+(\Phi_{\lambda,\mu})'(u)[v-u]\geq 0
\qquad\text{for every }v\in K_B.
\]
Notice that the constraint defining \(K_B\) only
guarantees \(\|\nabla u\|_{L^\infty(\Omega)}\leq 1\),
whereas solutions of the boundary value problem are required
to be \emph{strictly spacelike}, namely to satisfy
\(\|\nabla u\|_{L^\infty(\Omega)}<1\).
Accordingly, the correspondence between critical
points and solutions is far from being obvious.
The key observation to handle such an issue is that, if \(u\) is a critical
point of \(I_{\lambda,\mu}^B\), setting
\[
h_u(x):=\lambda u(x)+\mu \widetilde h(x,u(x))
\in L^\infty(\Omega),
\]
the above variational inequality means precisely
that \(u\) minimizes the convex functional
\[
v\mapsto\Psi_\mathcal B(v)-\int_\Omega h_u(x)v(x)\,dx
\]
over \(K_\mathcal B\). We are thus led to the auxiliary
problem
\[
\begin{cases}
-\displaystyle\operatorname{div}
\left(\dfrac{\nabla v}{\sqrt{1-|\nabla v|^2}}\right)
=h_u(x) & \text{in }\Omega,\\
\mathcal{B}v=0 & \text{on }\partial\Omega,
\end{cases}
\]
where, in the Neumann case, the compatibility
condition \(\int_\Omega h_u\,dx=0\) follows
directly from the variational inequality.
The required regularity and strict spacelikeness
results for the minimizers of this auxiliary problem
go back to Bartnik and Simon \cite{BaSi82} in the
Dirichlet case, while the corresponding Neumann
result has only recently been established by
Maniscalco \cite{Ma26} for convex domains. Building on these results, the arguments in
\cite[Theorem 2.1]{BeJeMa14} for the Dirichlet case and \cite[Theorem~2.10]{Ziegele} for the Neumann case can be adapted to the present
setting, yielding the following correspondence
between critical points of $I_{\lambda,\mu}^{\mathcal B}$ and strictly spacelike
solutions (in the weak sense at first, and then in the strong sense by standard elliptic regularity) of the associated nonlinear boundary value problem.

\begin{proposition}\label{prop:critical-points-solutions}
Let $\mathcal B\in\{\mathcal D,\mathcal N\}$ and assume
$(H_{\mathcal B})$, as well as $(h_1)$ and $(h_2)$. Moreover, let $\lambda,\mu\geq0$. Then, 
$u\in K_{\mathcal B}$ is a critical point of $I_{\lambda,\mu}^{\mathcal B}$ in the sense of Definition \ref{def-critical} if and only if $u$ is a strictly spacelike solution of
\begin{equation}
\label{eq:modified}
\begin{cases}
\displaystyle
-\operatorname{div}\left(
\frac{\nabla u}{\sqrt{1-|\nabla u|^2}}
\right)
=\lambda u+\mu \widetilde h(x,u)
& \text{in }\Omega,\\[1ex]
\mathcal B u=0
& \text{on }\partial\Omega.
\end{cases}
\end{equation}
In particular, every critical point of $I_{\lambda,\mu}^{\mathcal B}$ satisfies
$\|\nabla u\|_{L^\infty(\Omega)}<1.$
\end{proposition}

We next show that every solution of the truncated problem \eqref{eq:modified}
is also a solution of the original one.
For both boundary conditions, the key point is
an a priori bound on the solutions, obtained
for slightly different reasons.

\begin{proposition}\label{prop-truncation}
Let $\mathcal B\in\{\mathcal D,\mathcal N\}$ and assume
$(H_{\mathcal B})$, as well as $(h_1)$ and $(h_2)$.
Moreover, let $\lambda,\mu\geq0$.
Then every strictly spacelike solution of the
truncated problem \eqref{eq:modified} is also a solution
of the original problem \eqref{eq:problem},
with the same energy.

Moreover, for $\mathcal{B} = \mathcal{N}$ and
$(\lambda,\mu)\neq(0,0)$, every nontrivial
solution of \eqref{eq:modified} is sign-changing.
\end{proposition}

\begin{proof}
If $(\lambda,\mu)=(0,0)$, problems \eqref{eq:modified} and \eqref{eq:problem}
coincide, so there is nothing to prove.
Accordingly, we can assume from now on $(\lambda,\mu)\neq(0,0)$.

Let $u$ be a strictly spacelike solution of
the truncated problem \eqref{eq:modified}.
In particular, $u\in K_{\mathcal B}$ and
\[
|u(x)-u(y)|\leq |x-y|
\qquad\text{for every }x,y\in\overline\Omega.
\]
Indeed, in the Dirichlet case this follows
by extending $u$ by zero outside $\Omega$,
whereas in the Neumann case it follows from
the convexity of $\Omega$.

If $\mathcal B=\mathcal D$, since $u=0$
on $\partial\Omega$, we immediately obtain
\[
|u(x)|\leq\operatorname{dist}(x,\partial\Omega)
\leq\frac d2
\qquad\text{for every }x\in\Omega.
\]
On the other hand, for $\mathcal B=\mathcal N$, integrating the equation over $\Omega$ yields
\[
\int_\Omega
\bigl(\lambda u+\mu\widetilde h(x,u)\bigr)\,dx=0.
\]
Since $\widetilde h$ satisfies $(h_2)$,
every nontrivial solution must change sign.
Hence, by continuity, there exists
$y\in\Omega$ such that $u(y)=0$.
Consequently,
\[
|u(x)|\leq |x-y|\leq d
\qquad\text{for every }x\in\Omega.
\]
The same bound clearly holds if $u\equiv0$.

We have thus proved that, for both boundary
conditions,
\[
\|u\|_{L^\infty(\Omega)}\leq d.
\]
By the definition of the truncation,
\[
\widetilde h(x,u(x))=h(x,u(x)),
\qquad
\widetilde H(x,u(x))=H(x,u(x))
\quad\text{for a.e. }x\in\Omega.
\]
Therefore, $u$ solves the original problem \eqref{eq:problem}
with the same energy.
\end{proof}

\begin{remark}
Under the present assumptions, the converse
of Proposition \ref{prop-truncation} also holds,
as follows by applying the same argument to
solutions of the original problem.

We also observe that, since our variational
arguments are entirely carried out for the
truncated problem, assumption $(h_2)$ in
Theorems \ref{th-spectral} and \ref{th-mularge}
can actually be weakened. More precisely,
it suffices to require oddness and the strict
sign condition only for $0<s\leq d/2$ in
the Dirichlet case, and for $0<s\leq d$
in the Neumann case. Indeed, truncating the
nonlinearity at the corresponding threshold
preserves these properties, while the a priori
bounds established in the preceding proof
ensure that the solutions thus obtained also
solve the original problem.

Under these weaker assumptions, however,
the converse implication in Proposition
\ref{prop-truncation} need not hold in the
Neumann case.
\end{remark}

As a last result of this section, we show that, for either boundary condition,
the functional is bounded from below on a
suitable subspace of finite codimension. We also verify the weak Palais--Smale condition,
which will allow us to apply the abstract
multiplicity results of Section \ref{sec-abstract}.
\begin{lemma}\label{lem:compactness}
Let $\mathcal B\in\{\mathcal D,\mathcal N\}$ and assume
$(H_{\mathcal B})$, as well as $(h_1)$ and $(h_2)$. Moreover, let $\lambda,\mu\geq0$. 
Define
\[
Z_{\mathcal D}:=X,
\qquad
Z_{\mathcal N}
:=
\left\{
u\in X:\int_\Omega u\,dx=0
\right\}.
\]
Then
\[
\inf_{u\in Z_{\mathcal B}}
I_{\lambda,\mu}^{\mathcal B}(u)>-\infty.
\]
Moreover, $I_{\lambda,\mu}^{\mathcal B}$ satisfies the
weak Palais--Smale condition with respect to
$Y=C(\overline\Omega)$, provided that, in the Neumann
case, $(\lambda,\mu)\neq(0,0)$.
\end{lemma}

\begin{proof}
Let
\[
C_{\mathcal D}:=\frac d2,
\qquad
C_{\mathcal N}:=d.
\]
Arguing as in the proof of Proposition
\ref{prop-truncation}, the zero boundary
condition in the Dirichlet case, and the
zero-mean condition together with the
convexity of $\Omega$ in the Neumann case,
yield
\[
\|u\|_{L^\infty(\Omega)}\leq C_{\mathcal B}
\qquad
\text{for every }u\in K_{\mathcal B}\cap Z_{\mathcal B}.
\]
Since $\widetilde H(x,s)=H(x,s)$ for
$|s|\le d$, assumption $(h_1)$ yields
\[
I_{\lambda,\mu}^{\mathcal B}(u)
\geq
-\frac{\lambda}{2}|\Omega|C_{\mathcal B}^2
-\mu C_{\mathcal B}
\|\alpha_{C_{\mathcal B}}\|_{L^1(\Omega)}
\]
for every $u\in K_{\mathcal B}\cap Z_{\mathcal B}$.
Outside $K_{\mathcal B}$ the functional equals $+\infty$,
which proves the asserted lower bound.

Let $(u_n)\subset K_{\mathcal B}$ satisfy
$I_{\lambda,\mu}^{\mathcal B}(u_n)\to c\in\mathbb R$ and
\begin{equation}\label{eq:unified-wps-sequence}
\Psi_{\mathcal B}(v)-\Psi_{\mathcal B}(u_n)
+(\Phi_{\lambda,\mu})'(u_n)[v-u_n]
\geq-\varepsilon_n\|v-u_n\|_X
\qquad\text{for every }v\in X,
\end{equation}
where $\varepsilon_n\to0^+$.
We first show that $(u_n)$ is bounded in $X$.
For $\mathcal B=\mathcal D$, this follows directly
from the boundedness of $K_{\mathcal D}$. Suppose now that $\mathcal B=\mathcal N$: taking $v=u_n+1$ and $v=u_n-1$ in
\eqref{eq:unified-wps-sequence}, and using that $\Psi_{\mathcal B}(u+c) = \Psi_{\mathcal B}(u)$ for any $c \in \R$, we obtain
\begin{equation}\label{eq:approximate-compatibility}
\left|
\int_\Omega
\bigl(\lambda u_n+\mu\widetilde h(x,u_n)\bigr)\,dx
\right|
\leq\varepsilon_n.
\end{equation}
Now, write
$$u_n = v_n + m_n,$$
where $$m_n = \frac{1}{|\Omega|} \int_\Omega u_n \, dx,$$
and $v_n \in K_\mathcal N$ is such that $\int_\Omega v_n \, dx = 0$. By the convexity of $\Omega$ and the gradient constraint, one can deduce that
\[
|v_n| \leq d.
\]
If $(u_n)$ were unbounded in $L^\infty(\Omega)$,
then $(m_n)$ would be unbounded. Passing to a subsequence,
we would have either $m_n\to+\infty$ or $m_n\to-\infty$.
Since $\|v_n\|_{L^\infty(\Omega)}\leq d$, it would follow
that $u_n\to+\infty$ or $u_n\to-\infty$, respectively,
uniformly on $\overline\Omega$. In the first case, for all sufficiently large $n$,
\[
\widetilde h(x,u_n(x))=h(x,d)
\quad\text{a.e. in }\Omega.
\]
If $\lambda>0$, the integral in
\eqref{eq:approximate-compatibility} tends to $+\infty$.
If $\lambda=0$, then $\mu>0$ and that integral equals
\[
\mu\int_\Omega h(x,d)\,dx>0,
\]
again contradicting
\eqref{eq:approximate-compatibility}.
The case of uniform convergence to $-\infty$ is
analogous, using $h(x,-d)<0$.
Thus $(u_n)$ is bounded in $L^\infty(\Omega)$ and,
by the gradient constraint, also in $X$.

In either case, by the Arzel\`a--Ascoli theorem
and weak-* compactness, passing to a subsequence,
there exists $u\in K_{\mathcal B}$ such that
\[
u_n\to u\quad\text{in }C(\overline\Omega),
\qquad
\nabla u_n\rightharpoonup^\ast\nabla u
\quad\text{in }L^\infty(\Omega;\mathbb R^N).
\]
By the weak-* convergence of the gradients and
the convexity of the integrand (cf. \cite[Lemma 4 (iii)]{ABT})  we obtain
\[
\Psi_{\mathcal B}(u)
\leq\liminf_{n\to\infty}\Psi_{\mathcal B}(u_n).
\]
On the other hand, taking $v=u$ in
\eqref{eq:unified-wps-sequence} yields
\[
\Psi_{\mathcal B}(u_n)
\leq
\Psi_{\mathcal B}(u)
+(\Phi_{\lambda,\mu})'(u_n)[u-u_n]
+\varepsilon_n\|u-u_n\|_X.
\]
The last term tends to zero because $(u_n)$ is bounded
in $X$. Choose $R>0$ such that
\[
\|u\|_{L^\infty(\Omega)}\leq R,
\qquad
\sup_n\|u_n\|_{L^\infty(\Omega)}\leq R.
\]
By $(h_1)$ and the definition of the truncation, there
exists $a_R\in L^\infty(\Omega)$ such that
\[
|\widetilde h(x,s)|\leq a_R(x)
\qquad\text{for a.e. }x\in\Omega,\quad |s|\leq R.
\]
Consequently,
\[
\left|
(\Phi_{\lambda,\mu})'(u_n)[u-u_n]
\right|
\leq
\left(
\lambda R|\Omega|+\mu\|a_R\|_{L^1(\Omega)}
\right)
\|u-u_n\|_{L^\infty(\Omega)}
\longrightarrow0.
\]
It follows that
\[
\Psi_{\mathcal B}(u_n)\longrightarrow\Psi_{\mathcal B}(u).
\]
Uniform convergence and $(h_1)$ also imply
\[
\Phi_{\lambda,\mu}(u_n)
\longrightarrow\Phi_{\lambda,\mu}(u),
\]
and hence $I_{\lambda,\mu}^{\mathcal B}(u)=c$. Finally, for every fixed $v\in K_{\mathcal B}$,
dominated convergence gives
\[
(\Phi_{\lambda,\mu})'(u_n)[v-u_n]
\longrightarrow
(\Phi_{\lambda,\mu})'(u)[v-u].
\]
Since $\varepsilon_n\|v-u_n\|_X\to0$, passing to the limit
in \eqref{eq:unified-wps-sequence} yields
\[
\Psi_{\mathcal B}(v)-\Psi_{\mathcal B}(u)
+(\Phi_{\lambda,\mu})'(u)[v-u]\geq0.
\]
For $v\notin K_{\mathcal B}$, the inequality is automatic.
Thus $u$ is a critical point at level $c$, proving
the weak Palais--Smale condition.
\end{proof}

\begin{remark}
We point out that the truncation of the nonlinearity
is needed in the proof of the weak Palais--Smale
condition only in the Neumann case with $\lambda=0$.
Indeed, assumption $(h_2)$ alone does not prevent
$h(x,s)$ from approaching zero as $|s|\to\infty$,
which may lead to a loss of compactness along
the constant direction: the truncation rules out
this possibility. For simplicity, we nevertheless work with the
truncated nonlinearity in all cases, thus keeping
a unified variational framework throughout the paper.
\end{remark}

\section{Proof of Theorem~\ref{th-spectral}} \label{sec-3}

Fix $\mathcal B\in\{\mathcal D,\mathcal N\}$,
$\lambda>\lambda_k^{\mathcal B}$ with $k \geq 1$ and $\mu\geq0$.
Let us consider $Z_{\mathcal B}$ as in
Lemma~\ref{lem:compactness}. From that lemma, we know that
$I_{\lambda,\mu}^{\mathcal B}$ satisfies $\mathrm{(wPS)}$ and
\[
\inf_{u\in Z_{\mathcal B}}
I_{\lambda,\mu}^{\mathcal B}(u)>-\infty.
\]
Moreover,
\[
\operatorname{codim}_X Z_{\mathcal D}=0,
\qquad
\operatorname{codim}_X Z_{\mathcal N}=1.
\]
Now, let us define
\[
E_k^{\mathcal D}
:=
\operatorname{span}
\{\varphi_1^{\mathcal D},\ldots,\varphi_k^{\mathcal D}\},
\qquad
E_k^{\mathcal N}
:=
\operatorname{span}
\{\varphi_0^{\mathcal N},
  \varphi_1^{\mathcal N},\ldots,\varphi_k^{\mathcal N}\},
\]
where $\{\varphi_i^{\mathcal D}\}_{i=1}^{\infty}$ and
$\{\varphi_i^{\mathcal N}\}_{i=0}^{\infty}$ are
$L^2(\Omega)$-orthonormal eigenfunctions of the Dirichlet
and Neumann Laplacian, respectively, with
$\varphi_0^{\mathcal N}=|\Omega|^{-1/2}$.
By the regularity of the eigenfunctions,
$E_k^{\mathcal B}\subset X$; moreover, by the spectral decomposition of the Laplacian,
\[
\int_\Omega |\nabla u|^2\,dx
\leq \lambda_k^{\mathcal B}
\int_\Omega u^2\,dx
\qquad\text{for every }u\in E_k^{\mathcal B}.
\]
Finally, since $E_k^{\mathcal B}$ is finite-dimensional,
there exists a constant $C_k^{\mathcal B}>0$ such that
\[
\|\nabla u\|_{L^\infty(\Omega)}
\leq C_k^{\mathcal B}\|u\|_{L^\infty(\Omega)}
\qquad\text{for every }u\in E_k^{\mathcal B}.
\]
Let
\[
S_k^{\mathcal B}
:=
\left\{
u\in E_k^{\mathcal B}:
\|u\|_{L^\infty(\Omega)}=\rho
\right\},
\]
for some $\rho>0$ to be chosen later.
These sets are compact and symmetric in $X$, and
\[
S_k^{\mathcal D}\cong\mathbb S^{k-1},
\qquad
S_k^{\mathcal N}\cong\mathbb S^k,
\]
where the homeomorphisms can be chosen odd.

Now, fix $\delta\in(0,1)$. Note that
\[
1-\sqrt{1-|\nabla u|^2}
=
\frac{|\nabla u|^2}
{1+\sqrt{1-|\nabla u|^2}}.
\]
Hence, if $|\nabla u|\leq\delta$, then
\[
1-\sqrt{1-|\nabla u|^2}
\leq
\frac{1}{1+\sqrt{1-\delta^2}}|\nabla u|^2.
\]
Set
\[
C_\delta:=\frac{1}{1+\sqrt{1-\delta^2}}.
\]
If $u\in E_k^{\mathcal B}$ and
$\|\nabla u\|_{L^\infty(\Omega)}\leq\delta$, then
$u\in K_{\mathcal B}$. Using the fact that
$\widetilde H(x,s)\geq0$ for a.e. $x$ and every $s$ together with the spectral inequality
on $E_k^{\mathcal B}$, we obtain
\begin{align*}
I_{\lambda,\mu}^{\mathcal B}(u)
&=
\int_\Omega
\left(1-\sqrt{1-|\nabla u|^2}\right)\,dx
-\frac{\lambda}{2}\int_\Omega u^2\,dx
-\mu\int_\Omega \widetilde H(x,u)\,dx\\
&\leq
C_\delta\int_\Omega|\nabla u|^2\,dx
-\frac{\lambda}{2}\int_\Omega u^2\,dx\\
&\leq
\left(
C_\delta\lambda_k^{\mathcal B}-\frac{\lambda}{2}
\right)\|u\|_{L^2(\Omega)}^2.
\end{align*}
Since $C_\delta\to\frac12$ as $\delta\to0$ and
$\lambda>\lambda_k^{\mathcal B}$, we can choose
$\delta$ sufficiently small so that
\[
C_\delta\lambda_k^{\mathcal B}-\frac{\lambda}{2}<0.
\]
We then choose $\rho \in (0,\frac{\delta}{C_k^{\mathcal B}})$, so that
\[
u\in E_k^{\mathcal B},
\quad
\|u\|_{L^\infty(\Omega)}=\rho
\quad\Longrightarrow\quad
\|\nabla u\|_{L^\infty(\Omega)}<\delta.
\]
With this choice, $S_k^{\mathcal B}\subset K_{\mathcal B}$.
Moreover, by compactness of $S_k^{\mathcal B}$ and
since $0\notin S_k^{\mathcal B}$,
\[
\beta_k^{\mathcal B}
:=
\min_{u\in S_k^{\mathcal B}}
\|u\|_{L^2(\Omega)}^2>0.
\]
Consequently,
\[
\sup_{u\in S_k^{\mathcal B}}
I_{\lambda,\mu}^{\mathcal B}(u)
\leq
\left(
C_\delta\lambda_k^{\mathcal B}-\frac{\lambda}{2}
\right)\beta_k^{\mathcal B}
<0=I_{\lambda,\mu}^{\mathcal B}(0).
\]

We may therefore apply Theorem~\ref{thm:ABT-negative}, with
\[
Z=Z_{\mathcal B},
\qquad
A=S_k^{\mathcal B},
\qquad
(\ell,m)=
\begin{cases}
(1,k), & \mathcal B=\mathcal D,\\
(2,k+1), & \mathcal B=\mathcal N.
\end{cases}
\]
In both cases, we obtain at least
\[
m-\ell+1
=
\dim E_k^{\mathcal B}
-\operatorname{codim}_X Z_{\mathcal B}
=k
\]
distinct pairs of nonzero critical points of
$I_{\lambda,\mu}^{\mathcal B}$ with negative energy. 

By Propositions \ref{prop:critical-points-solutions} and
\ref{prop-truncation}, these critical points
yield distinct pairs of strictly spacelike
solutions of problem \eqref{eq:problem}. In the Neumann case, all these
solutions are sign-changing.
In the Dirichlet case, moreover, since 
$I_{\lambda,\mu}^\mathcal{D}(u) = I_{\lambda,\mu}^\mathcal{D}(|u|)$, the pair $(u,-u)$ at level $c_1 = \inf_X I_{\lambda,\mu}^\mathcal{D}$ can be chosen so that $u$ has constant sign. The proof is thus concluded.

\section{Proof of Theorem \ref{th-mularge}} \label{sec-4}

Fix $\mathcal B\in\{\mathcal D,\mathcal N\}$, $\lambda\in
[\lambda_k^{\mathcal B},\lambda_{k+1}^{\mathcal B})$ for some integer $k\geq0$, and an integer
$j\geq\max\{1,k\}$.
We first obtain $j$ pairs of solutions with negative
energy and then, if $j>k$, other $j-k$ pairs with
positive energy.

Consider the spaces
\[
E_j^{\mathcal D}
:=
\operatorname{span}
\{\varphi_1^{\mathcal D},\ldots,\varphi_j^{\mathcal D}\},
\qquad
E_j^{\mathcal N}
:=
\operatorname{span}
\{\varphi_0^{\mathcal N},
  \varphi_1^{\mathcal N},\ldots,\varphi_j^{\mathcal N}\},
\]
where the eigenfunctions are those introduced in
the proof of Theorem \ref{th-spectral}. In particular,
$E_j^{\mathcal B}\subset X$, and there exists
$C_j^{\mathcal B}>0$ such that
\[
\|\nabla u\|_{L^\infty(\Omega)}
\leq C_j^{\mathcal B}\|u\|_{L^\infty(\Omega)}
\qquad\text{for every }u\in E_j^{\mathcal B}.
\]
Choose now $r \in (0,1/C_j^{\mathcal{B}})$
and define
\[
S_j^{\mathcal B}
:=
\left\{
u\in E_j^{\mathcal B}:
\|u\|_{L^\infty(\Omega)}=r
\right\}.
\]
Then $S_j^{\mathcal B}\subset K_{\mathcal B}$.
Moreover, $S_j^{\mathcal B}$ is compact and symmetric
in $X$, and is oddly homeomorphic to
$\mathbb S^{j-1}$ if $\mathcal B=\mathcal D$, and to
$\mathbb S^j$ if $\mathcal B=\mathcal N$. 
By $(h_1)$, the map
\[
u\longmapsto\int_\Omega \widetilde H(x,u)\,dx
\]
is continuous on this set.
Furthermore, $(h_2)$ implies that $  \widetilde H(x,s)>0$ for
$s\neq0$ and a.e. $x\in\Omega$.
Since every element of $S_j^{\mathcal B}$ is nonzero,
compactness yields
\[
M_j^{\mathcal B}
:=
\min_{u\in S_j^{\mathcal B}}
\int_\Omega  \widetilde H(x,u)\,dx>0.
\]
Also, set
\[
b_j^{\mathcal B}(\lambda)
:=
\max_{u\in S_j^{\mathcal B}}
I_{\lambda,0}^{\mathcal B}(u)\in\mathbb R.
\]
Choosing
\[
\mu_{j,\mathcal B}^{\star}(\lambda)
:=
1+
\frac{\max\{0,b_j^{\mathcal B}(\lambda)\}}
     {M_j^{\mathcal B}},
\]
for every $\mu>\mu_{j,\mathcal B}^{\star}(\lambda)$,
we then have
\begin{align*}
\sup_{u\in S_j^{\mathcal B}}
I_{\lambda,\mu}^{\mathcal B}(u)
&=
\sup_{u\in S_j^{\mathcal B}}
\left(
I_{\lambda,0}^{\mathcal B}(u)
-\mu\int_\Omega  \widetilde H(x,u)\,dx
\right)\\
&\leq b_j^{\mathcal B}(\lambda)
-\mu M_j^{\mathcal B}<0.
\end{align*}
Fix such a $\mu$.
By Lemma~\ref{lem:compactness},
$I_{\lambda,\mu}^{\mathcal B}$ satisfies
$\mathrm{(wPS)}$ and
\[
\inf_{u\in Z_{\mathcal B}}
I_{\lambda,\mu}^{\mathcal B}(u)>-\infty.
\]
Applying Theorem~\ref{thm:ABT-negative} with
\[
Z=Z_{\mathcal B},
\qquad
A=S_j^{\mathcal B},
\qquad
(\ell,m)=
\begin{cases}
(1,j), & \mathcal B=\mathcal D,\\
(2,j+1), & \mathcal B=\mathcal N,
\end{cases}
\]
we obtain $j$ distinct pairs of nonzero critical points
with negative energy.

Suppose now that $j>k$.
To obtain critical points with positive energy,
we apply Theorem~\ref{thm:ABT-positive}.
With the conventions
\[
E_0^{\mathcal D}:=\{0\},
\qquad
E_0^{\mathcal N}
:=\operatorname{span}\{\varphi_0^{\mathcal N}\},
\]
define
\[
W_k^{\mathcal B}
:=
\left\{
u\in X:
\int_\Omega uv\,dx=0
\quad\text{for every }v\in E_k^{\mathcal B}
\right\}.
\]
This is a closed linear subspace of $X$, and
\[
\operatorname{codim}_X W_k^{\mathcal D}=k,
\qquad
\operatorname{codim}_X W_k^{\mathcal N}=k+1.
\]
For every $u\in W_k^{\mathcal B}\cap K_{\mathcal B}$,
the spectral decomposition of the Laplacian gives
\[
\int_\Omega|\nabla u|^2\,dx
\geq
\lambda_{k+1}^{\mathcal B}\int_\Omega u^2\,dx.
\]
Set $\sigma:=\lambda_{k+1}^{\mathcal B}-\lambda>0$. By $(h_3)$, there exists $r_\mu\in(0,d)$ such that
\[
0\leq \widetilde H(x,s)=H(x,s)
\leq\frac{\sigma}{4\mu}s^2 \quad \text{for a.e. $x\in\Omega$ and every $|s|\leq r_\mu$.}
\]
Using
\[
1-\sqrt{1-t^2}\geq\frac12t^2
\qquad\text{for }0\leq t\leq1,
\]
we obtain, for every
$u\in W_k^{\mathcal B}\cap K_{\mathcal B}$
with $\|u\|_{L^\infty(\Omega)}\leq r_\mu$,
\begin{align*}
I_{\lambda,\mu}^{\mathcal B}(u)
&\geq
\frac12\int_\Omega|\nabla u|^2\,dx
-\frac{\lambda}{2}\int_\Omega u^2\,dx
-\mu\int_\Omega \widetilde H(x,u)\,dx\\
&\geq
\frac{\sigma}{2}\|u\|_{L^2(\Omega)}^2
-\frac{\sigma}{4}\|u\|_{L^2(\Omega)}^2\\
&=
\frac{\sigma}{4}\|u\|_{L^2(\Omega)}^2.
\end{align*}
Choose $$0<\rho<\min \left\{r,r_\mu, \frac{\|\varphi_{k+1}^\mathcal B \|_{L^\infty(\Omega)}}{\| \nabla \varphi_{k+1}^\mathcal B \|_{L^\infty (\Omega)}}\right\},$$
and define
\[
\mathcal S_{k,\rho}^{\mathcal B}
:=
\left\{
u\in W_k^{\mathcal B}\cap K_{\mathcal B}:
\|u\|_{L^\infty(\Omega)}=\rho
\right\}.
\]
Note that $\mathcal S_{k,\rho}^{\mathcal B}$ is non-empty, since $u = \rho \frac{\varphi_{k+1}^\mathcal B}{\|\varphi_{k+1}^\mathcal B\|_{L^\infty (\Omega)}} \in\mathcal S_{k,\rho}^{\mathcal B}$. 
Since the functions in $S_{k,\rho}^{\mathcal B}$
are uniformly bounded and $1$-Lipschitz
(cf. the proof of Proposition \ref{prop-truncation}),
the Arzel\`a--Ascoli theorem yields that
$S_{k,\rho}^{\mathcal B}$ is relatively compact
in $C(\overline\Omega)$.
Moreover, $S_{k,\rho}^{\mathcal B}$ is closed
with respect to uniform convergence,
and hence compact.
Since $0\notin\mathcal S_{k,\rho}^{\mathcal B}$, we conclude that
\[
\beta_{k,\rho}^{\mathcal B}
:=
\min_{u\in\mathcal S_{k,\rho}^{\mathcal B}}
\|u\|_{L^2(\Omega)}^2>0.
\]
Consequently, setting
\[
\alpha:=
\frac{\sigma}{4}\beta_{k,\rho}^{\mathcal B}>0,
\]
we have
\[
I_{\lambda,\mu}^{\mathcal B}(u)\geq\alpha
\qquad
\text{for every }u\in W_k^{\mathcal B}
\text{ with }\|u\|_{L^\infty(\Omega)}=\rho.
\]
Indeed, the preceding estimate applies when
$u\in K_{\mathcal B}$, whereas the functional equals
$+\infty$ otherwise.

The hypotheses of Theorem~\ref{thm:ABT-positive}
are therefore satisfied with
\[
E=E_j^{\mathcal B},
\qquad
Z=W_k^{\mathcal B},
\qquad
(\ell,m)=
\begin{cases}
(k,j), & \mathcal B=\mathcal D,\\
(k+1,j+1), & \mathcal B=\mathcal N.
\end{cases}
\]
We obtain at least $m-\ell=j-k$ distinct pairs
of nonzero critical points with positive energy.
These are distinct from the $j$ pairs obtained at
negative energy, giving at least $2j-k$ pairs in total.
If $j=k$, the negative-energy construction already
gives the required $j=2j-k$ pairs.

By Propositions \ref{prop:critical-points-solutions} and
\ref{prop-truncation}, these critical points
yield distinct pairs of strictly spacelike
solutions of problem \eqref{eq:problem}. In the Neumann case, all these
solutions are sign-changing. In the Dirichlet case, moreover, a pair of
negative-energy solutions can be chosen to
have constant sign, by the same argument
used in the proof of Theorem \ref{th-spectral}.


\begin{thebibliography}{99}


\bibitem{ABT}
D.~Arcoya, C.~Bereanu, and P.~J.~Torres,
Lusternik--Schnirelman theory for the action integral of the Lorentz force equation,
\emph{Calc. Var. Partial Differential Equations}
\textbf{59} (2020), Paper No.~50, 32 pp.


\bibitem{BaSi82}
R.~Bartnik and L.~Simon,
Spacelike hypersurfaces with prescribed boundary values and mean curvature,
\emph{Comm. Math. Phys.}
\textbf{87} (1982/83), 131--152.

\bibitem{BeJeMa14}
C.~Bereanu, P.~Jebelean, and J.~Mawhin,
The Dirichlet problem with mean curvature operator in Minkowski space:
a variational approach,
\emph{Adv. Nonlinear Stud.}
\textbf{14} (2014), 315--326.

\bibitem{BeJeTo13}
C.~Bereanu, P.~Jebelean, and P.~J.~Torres,
Positive radial solutions for Dirichlet problems with mean curvature operators
in Minkowski space,
\emph{J. Funct. Anal.}
\textbf{264} (2013), 270--287.

\bibitem{BoFoMoSaTa18}
D.~Bonheure, J.~F\"oldes, E.~Moreira dos Santos, A.~Salda\~na, and H.~Tavares,
Paths to uniqueness of critical points and applications to partial differential equations,
\emph{Trans. Amer. Math. Soc.}
\textbf{370} (2018), no.~10, 7081--7127.


\bibitem{BoCoZi26}
A.~Boscaggin, F.~Colasuonno, and R.~Ziegele,
Positive and nodal solutions for the Minkowski mean curvature equation:
multiplicity and asymptotics,
preprint, arXiv:2607.15956 [math.AP] (2026).


\bibitem{CoObOmRi13b}
C.~Corsato, F.~Obersnel, P.~Omari, and S.~Rivetti,
 Positive solutions of the Dirichlet problem for the prescribed mean curvature equation
 in Minkowski space,
 \emph{J. Math. Anal. Appl.}
 \textbf{405} (2013), 227--239.

\bibitem{Da16}
G.~Dai,
Bifurcation and positive solutions for problem with mean curvature operator
in Minkowski space,
\emph{Calc. Var. Partial Differential Equations}
\textbf{55} (2016), no.~4, Paper No.~72.

\bibitem{GaLo26}
M.~Garz\'on and S.~L\'opez-Mart\'inez,
A local linking theorem for relativistic action functionals,
preprint, arXiv:2607.01845 [math.AP] (2026).


\bibitem{JS}
P.~Jebelean and C. \c{S}erban,
Fisher--Kolmogorov type perturbations of the mean curvature operator in Minkowski space,
\emph{Electron. J. Qual. Theory Differ. Equ.}
(2020), Paper No.~81, 1--12.

\bibitem{Ma26}
L.~Maniscalco,
Free boundary space-like graphs with prescribed mean curvature,
preprint, arXiv:2608.00887 [math.AP] (2026).


\bibitem{Mi21}
M.~Mih\u{a}ilescu,
The spectrum of the mean curvature operator,
\emph{Proc. Roy. Soc. Edinburgh Sect. A}
\textbf{151} (2021), no.~2, 451--463.


\bibitem{Sz86}
A.~Szulkin,
Minimax principles for lower semicontinuous functions and applications to nonlinear
boundary value problems,
\emph{Ann. Inst. H. Poincar\'e C Anal. Non Lin\'eaire}
\textbf{3} (1986), no.~2, 77--109.

\bibitem{Ziegele}
R.~Ziegele,
Positive solutions for an indefinite-weight Minkowski mean curvature Neumann problem:
multiplicity and asymptotic behaviour,
preprint, arXiv:2609.09894 [math.AP] (2026).

\end{thebibliography}
\end{document}